\documentclass[11pt]{article}%
\usepackage{amssymb}
\usepackage{amsmath}
\usepackage{amsfonts}
\usepackage{graphicx}

\newtheorem{theorem}{Theorem}
\newtheorem{definition}[theorem]{Definition}
\newtheorem{example}[theorem]{Example}

\newtheorem{proposition}[theorem]{Proposition}

\begin{document}

\title{Geometric structure for the tangent bundles of direct limit manifolds}

\author{A. Suri and P. Cabau}
\maketitle

\textbf{Abstract:} We equip the direct limit of tangent bundles of
paracompact finite dimensional manifolds with a structure of convenient vector
bundle with structural group $GL\left(  \infty,\mathbb{R}\right)
=\underrightarrow{\lim}GL\left(  \mathbb{R}^{n}\right)  $\textit{.}\\
\\
\textbf{R\'{e}sum\'{e} :} On munit la limite directe des fibr\'{e}s tangents
\`{a} des vari\'{e}t\'{e}s paracompactes de dimensions finies d'une structure
de fibr\'{e} vectoriel 'convenient' (au sens de Kriegel et Michor) de groupe
structural $GL\left(  \infty,\mathbb{R}\right)  =\underrightarrow{\lim}GL\left(  \mathbb{R}^{n}\right)  $\textit{.}
\\
\\
\textbf{MSC classification 2000:}    Primary 58B20; Secondary 58A05.

\section{Introduction}

G. Galanis proved in \cite{Gal} that the tangent bundle of a projective limit
of Banach manifolds can be equipped with a Fr\'{e}chet vector bundle structure
with structural group a topological subgroup of the general linear group of
the fiber type. Various problems were studied in this framework:
connections, ordinary differential equations, ... (\cite{ADGS}, \cite{AghSur}%
, ...).

Here we consider the situation for direct (or inductive) limit
of tangent bundles $TM_{i}$ where $M_{i}$ is a finite dimensional manifold: we
first have (Proposition \ref{P_StructureConvenientManifold}) that
$M=\underrightarrow{\lim}M_{i}$ can be endowed with a structure of convenient
manifold modelled on the convenient vector space $\mathbb{R}^{\infty}=$
$\underrightarrow{\lim}\mathbb{R}^{n}$ of finite sequences, equipped with the
finite topology (cf \cite{Han}). We then prove (Theorem
\ref{T_StructureConvenientVectorBundleOnTM}) that $TM$ can be endowed with a
convenient structure of vector bundle whose structural group is $GL\left(
\infty,\mathbb{R}\right)  =\underrightarrow{\lim}GL\left(  \mathbb{R}%
^{n}\right)$ (the group of invertible matrices of countable size, differing from the identity matrix at only finitely many places, first described by Milnor in \cite{Mil}). As an example we consider
the tangent bundle to $\mathbb{S}^{\infty}$. Other examples
can be found in the framework of manifolds for algebraic topology, such as
Grassmannians (\cite{KriMic}) or Lie groups (\cite{Glo1}, \cite{Glo2}).

The paper is organized as follows: We first recall the framework of convenient
calculus (part \ref{_ConvenientCalculus}). In part \ref{_DirectLimits}, we
review direct limit in different categories. We obtain the main result
(theorem \ref{T_StructureConvenientVectorBundleOnTM}) in the last part.

\section{\label{_ConvenientCalculus}Convenient calculus}

Classical differential calculus is perfectly adapted to finite dimensional or
even Banach manifolds (cf. \cite{Lan}).

On the other hand, convenient analysis, developed in \cite{KriMic}, provides a
satisfactory solution of the question how to do analysis on a large class of
locally convex spaces and in particular on strict inductive limits of Banach
manifolds or fiber bundles.

In order to endow some locally convex vector spaces (l.c.v.s.) $E$, which will
be assumed Hausdorff, with a differentiable structure we first use the notion
of smooth curves $c:\mathbb{R}\rightarrow E$, which poses no problems.

We denote the space $C^{\infty}\left(  \mathbb{R},E\right)  $ by $\mathcal{C}%
$; the set of continuous linear functionals is denoted by $E^{^{\prime}}$.

We then have the following characterization: a subset $B$ of $E$ is bounded
iff $l\left(  B\right)  $ is bounded for any $l\in E^{^{\prime}}$.

\begin{definition}
\label{D_MackeyCauchySequence}A sequence $\left(  x_{n}\right)  $ in $E$ is
called Mackey-Cauchy if there exists a bounded absolutely convex set $B$ and
for every $\varepsilon>0$ an integer $n_{\varepsilon}\in N$ s.t. $a_{n}%
-a_{m}\in\varepsilon B$ whenever $n>m>n_{\varepsilon}$
\end{definition}

\begin{definition}
\label{D_ConvenientVectorSpace}A locally convex vector space is said to be
$c^{\infty}$-complete or \textit{convenient }if one of the following
(equivalent) conditions is satisfied~:

\begin{enumerate}
\item if $c:\mathbb{R}\rightarrow E$ is a curve such that $l\circ
c:\mathbb{R}\rightarrow\mathbb{R}$ is smooth for all continuous linear
functionnal $l$, then $c$ is smooth.

\item Any Mackey-Cauchy sequence converges\footnote{This condition is
equivalent to:
\par
For every absolutely convex closed bounded set $B$ the linear span $E_{B}$ of
$B$ in $E$, equipped with the Minkowski functional $p_{B}\left(  v\right)
=\inf\left\{  \lambda>0:v\in\lambda.B\right\}  $, is complete.} (i.e. $E$ is
Mackey complete)

\item For any $c\in\mathcal{C}$ there exists $\gamma\in\mathcal{C}$ such that
$\gamma^{\prime}=c$.
\end{enumerate}
\end{definition}

The $c^{\infty}$--topology on a l.c.v.s. is the final topology with respect to
all smooth curves $\mathbb{R}\rightarrow E$ ; it is denoted by $c^{\infty}E$.
Its open sets will be called $c^{\infty}-$open.

Note that the $c^{\infty}$--topology is finer than the original topology. For
Fr\'{e}chet spaces, this topology coincides with the given locally convex topology.

In general, $c^{\infty}E$ is not a topological vector space.

The following theorem gives some constructions inheriting of $c^{\infty}$-completeness.

\begin{theorem}
\label{T_InheritanceCompleteness}The following constructions preserve
$c^{\infty}$-completeness: limits, direct sums, strict inductive limits of
sequences of closed embeddings.
\end{theorem}

\medskip The category $\mathbb{CON}$ of convenient vector spaces and
continuous linear maps forms a symmetric monoidal closed category.

Let $E$ and $F$ be two convenient spaces and let $U\subset E$ be a $c^{\infty
}-$open. A map $f:E\supset U\rightarrow F$ is said to be smooth if $f\circ
c\in C^{\infty}\left(  \mathbb{R},F\right)  $ for any $c\in C^{\infty}\left(
\mathbb{R},U\right)  .$ Moreover, the space $C^{\infty}\left(  U,F\right)  $
may be endowed with a structure of convenient vector space.

Let $L\left(  E,F\right)  $ be the space of all bounded linear mappings. We
can define the differential operator%
\begin{align*}
d  &  :C^{\infty}\left(  E,F\right)  \rightarrow C^{\infty}\left(  E,L\left(
E,F\right)  \right) \\
df\left(  x\right)  v  &  =\underset{t\rightarrow0}{\lim}\dfrac{f\left(
x+tv\right)  -f\left(  x\right)  }{t}%
\end{align*}
which is linear and bounded (and so smooth).

\section{\label{_DirectLimits}Direct (or inductive) limits}

\subsection{\bigskip\label{__DirectLimitCategory}Direct limit in a category}

The references are \cite{Bou} and \cite{Glo2}.

\begin{definition}
A direct sequence in a category $\mathbb{A}$ is a pair $\mathcal{S}=\left(
X_{i},\varepsilon_{ij}\right)  _{\left(  i,j\right)  \in\mathbb{N}^{2},\ i\leq
j}$ where $X_{i}$ is an object of $\mathbb{A}$ and each $\varepsilon
_{ij}:X_{i}\rightarrow X_{j}$ is a morphism, called bonding map, such
that:\newline\newline-- $\varepsilon_{ii}=\operatorname{Id}_{X_{i}}$%
\newline\newline-- $\varepsilon_{jk}\circ\varepsilon_{ij}=\varepsilon_{ik}$ if
$i\leq j\leq k$
\end{definition}

\begin{definition}
A cone over $\mathcal{S}$ is a pair $\left(  X,\varepsilon_{i}\right)
_{i\in\mathbb{N}}$ where $X$ is an object of $\mathbb{A}$ and $\varepsilon
_{i}:X_{i}\rightarrow X $ is a morphism of this category such that
\[
\varepsilon_{j}\circ\varepsilon_{ij}=\varepsilon_{i}\text{ if }i\leq j
\]

A cone $\left(  X,\varepsilon_{i}\right)  _{i\in\mathbb{N}}$ is a direct limit
cone over $\mathcal{S}$ in the category $\mathbb{A}$ if for every cone
$\left(  Y,\psi_{i}\right)  $ over $\mathcal{S}$ there exists a unique
morphism $\psi:X\rightarrow Y$ such that $\psi\circ\varepsilon_{i}=\psi_{i}$
for each $i$.\newline We then write $X=\underrightarrow{\lim}X_{i}$ and we
call $X$ the direct limit of $\mathcal{S}$.
\end{definition}

\subsection{\label{__DirectLimitSets}Direct limit of sets}

Let $\mathcal{S}=\left(  X_{i},\varepsilon_{ij}\right)  _{\left(  i,j\right)
\in\mathbb{N}^{2},\ i\leq j}$ be a direct sequence of sets.

\bigskip The \textit{direct sum }$\bigoplus\limits_{n\in\mathbb{N}}$ $X_{n}$
also called the \textit{coproduct} $\coprod\limits_{n\in\mathbb{N}}X_{n}$ is
the subspace of the cartesian product $\prod\limits_{n\in\mathbb{N}}X_{n}$
formed by all the points with only finetely many non-vanishing coordinates.

In this space we introduce the following binary relation (where $x\in X_{i}$
and $y\in Y_{j}$)%
\[
\left(  i,x\right)  \sim\left(  j,y\right)  \quad\Longleftrightarrow
\quad\left\vert
\begin{array}
[c]{c}%
y=\varepsilon_{ij}\left(  x\right)  \ \text{if }i\leq j\\
\text{or}\\
x=\varepsilon_{ji}\left(  y\right)  \text{ if }i\geq j
\end{array}
\right.
\]

which is an equivalence relation$.$

Then the set $X=\coprod\limits_{n\in\mathbb{N}}X_{n}\ /\sim$ together with the
maps
\[%
\begin{array}
[c]{cccc}%
\varepsilon_{i}: & X_{i} & \longrightarrow & X\\
& x & \mapsto & \widetilde{\left(  i,x\right)  }%
\end{array}
\]

where $\widetilde{\left(  i,x\right)  }$ is the equivalence class of $\left(
i,x\right)  $, is the \textit{direct limit} of $\mathcal{S}$ in the category
$\mathbb{SET}$.

We have $X=\bigcup\limits_{i\in N}\varepsilon_{i}\left(  X_{i}\right)  $. If
each $\varepsilon_{ij}$ is injective then so is $\varepsilon_{i}$.
$\mathcal{S}$ is then equivalent to the sequence of the subsets $\varepsilon
_{i}\left(  X_{i}\right)  \subset X$ with the inclusion maps.

\subsection{\label{__DirectLimitTopologicalSpaces}Direct limit of topological
spaces}

Let $\mathcal{S}=\left(  X_{i},\varepsilon_{ij}\right)  _{\left(  i,j\right)
\in\mathbb{N}^{2},\ i\leq j}$ be a direct sequence of topological spaces where
the bonding maps are continuous.

We then endow $X$ with the direct sum topology, i.e. is the final topology
with respect to the family $\left(  \varepsilon_{i}\right)  _{i\in\mathbb{N}}$
which is the finest topology for which the maps $\varepsilon_{i}$ are
continuous. Then $U\subset X$ is open if and only if $\left(  \varepsilon
_{i}\right)  ^{-1}\left(  U\right)  $ is open in $X_{i}$ for each $i$.

If the bonding maps are topological embeddings we call $\mathcal{S}$
\textit{strict direct limi}t. For any $i\in\mathbb{N}$, $\varepsilon_{i}$ is
then a topological embedding.

\subsection{\label{__exampleSequenceFinitelyManyNonZerosTerms}Fundamental
example of $\mathbb{R}^{\infty}$}

The space $\mathbb{R}^{\infty}$ also denoted by $\mathbb{R}^{\left(
\mathbb{N}\right)  }$ of all finite sequences is the direct limit of $\left(
\mathbb{R}^{i},\varepsilon_{ij}\right)  _{\left(  i,j\right)  \in
\mathbb{N}^{2},\ i\leq j}$ where $\varepsilon_{ij}:\left(  x_{1},\dots
,x_{i}\right)  \mapsto\left(  x_{1},\dots,x_{i},0,\dots,0\right)  .$

It is a convenient vector space (\cite{KriMic}, 47.1).

\subsection{Direct limit of finite dimensional manifolds}

Let $\mathcal{M}=\left(  M_{i},\phi_{ij}\right)  _{\ i\leq j}$ be a direct
sequence of paracompact finite dimensional smooth real manifolds where the
bonding maps $\phi_{ij}:M_{i}\longrightarrow M_{j}$ are injective smooth
immersions and $\underset{i\in\mathbb{N}}{\sup}\left\{  \dim_{\mathbb{R}}%
M_{i}\right\}  =\infty$. Adapting a result of Gl\"{o}ckner (\cite{Glo2},
Theorem 3.1) to the convenient framework (using Proposition 3.6) we have:

\begin{theorem}
There exists a uniquely determined $c^{\infty}-$manifold structure on the
direct limit $M$ of $\mathcal{M}$ modelled on the convenient vector space
$\mathbb{R}^{\infty}$.
\end{theorem}

\begin{example}
The sphere $\mathbb{S}^{\infty}$ (\cite{KriMic}, 47.2).-- The convenient
vector space $\mathbb{R}^{\infty}$ is equipped with the weak inner product
given by the finite sum $\left\langle x,y\right\rangle =\sum\limits_{i}%
x_{i}y_{i}$ and is bilinear and bounded, therefore smooth. The topological
inductive limit of $\mathbb{S}^{1}\subset\mathbb{S}^{2}\subset\cdots$ is the
closed subset $\mathbb{S}^{\infty}=\left\{  x\in\mathbb{R}^{\infty
}:\left\langle x,x\right\rangle =1\right\}  $ of $\mathbb{R}^{\infty}$.
\newline\newline Choose $a\in\mathbb{S}^{\infty}$. We can define the
stereographic atlas corresponding to the equivalence class of the two charts
$\left\{  \left(  U_{+},u_{+}\right)  ,\left(  U_{-},u_{-}\right)  \right\}  $
where $U_{+}=\mathbb{S}^{\infty}\backslash\left\{  a\right\}  $ (resp.
$U_{-}=\mathbb{S}^{\infty}\backslash\left\{  -a\right\}  $) and $%
\begin{array}
[c]{cccc}%
u_{+}: & U_{+} & \longrightarrow & \left\{  a\right\}  ^{\perp}\\
& x & \mapsto & \dfrac{x-\left\langle x,a\right\rangle a}{1-\left\langle
x,a\right\rangle }%
\end{array}
$ (resp. $%
\begin{array}
[c]{cccc}%
u_{-}: & U_{-} & \longrightarrow & \left\{  a\right\}  ^{\perp}\\
& x & \mapsto & \dfrac{x-\left\langle x,a\right\rangle a}{1+\left\langle
x,a\right\rangle }%
\end{array}
$). Then $\mathbb{S}^{\infty}$ is a convenient manifold modelled on
$\mathbb{R}^{\infty}$.
\end{example}

\section{\label{_TangentBundleDirectLimitManifolds} Tangent bundle of direct limit of manifolds}

\subsection{Structure of manifold on direct limit of tangent bundles}

Let $p\geq4$ and $\{M_{i},\phi_{ij}\}_{i\leq j}$ be a direct sequence of
$C^{p}$ paracompact finite dimensional manifolds for which the connecting
morphisms are $C^{p}$ embeddings with closed image. Without loss of generality
(cf. \ref{__DirectLimitSets}) we may assume that $M_{1}\subseteq
M_{1}\subseteq\dots\subseteq M$ where $\{M,\phi_{i}\}$ is the direct limit of
$\{M_{i},\phi_{ij}\}_{i\leq j}$ in the category of topological spaces and the
maps $\phi_{i}:M_{i}\longrightarrow M$ are inclusions \cite{Glo2}. Suppose
that $\dim M_{i}=d_{i}$ and consider for $i\leq j$,%
\[%
\begin{array}
[c]{cccc}%
\lambda_{ij}: & \mathbb{R}^{d_{i}} & \longrightarrow & \mathbb{R}^{d_{j}}\\
& \left(  x_{1},\dots,x_{d_{i}}\right)  & \mapsto & \left(  x_{1}%
,\dots,x_{d_{i}},0,\dots,0\right)
\end{array}
\]
For $x\in M$ there exists $n\in\mathbb{N}$ such that $x=\phi_{n}(x)$. Using
tubular neighborhoods Gl$\ddot{\text{o}}$ckner proved that there exists an
open neighborhood $O_{x}$ of $x$ in $M$ and a sequence of $C^{p-2}$
diffeomorphisms $\{h_{i}^{(x)}:\mathbb{R}^{d_{i}}\longrightarrow
U_{i}\}_{i\geq n}$ (inverse of chart mappings) where $U_{i}={\phi_{i}}%
^{-1}(O_{x})$. Moreover for $j\geq i\geq n$ the compatibility condition
\begin{equation}
h_{j}^{(x)}\circ\lambda_{ij}=\phi_{ij}|_{U_{i}}\circ h_{i}^{(x)}
\label{compatibility condition of hi}%
\end{equation}
holds true (\cite{Glo1}, Lemma 4.1).

Our first aim is to introduce appropriate connecting morphisms, say
$\{\Phi_{ij}\}_{i\leq j}$, such that $\{TM_{i},\Phi_{ij}\}$ form a direct
system of manifolds in the sense of Gl$\ddot{\text{o}}$ckner.

For $i\leq j$ define
\begin{align*}
\Phi_{ij}:TM_{i}  &  \longrightarrow TM_{j}\\
{[\alpha_{i},x_{i}]_{i}}  &  \longmapsto\lbrack\phi_{ij}\circ\alpha_{i}%
,\phi_{ij}(x_{i})]_{j}%
\end{align*}
where the bracket $\left[  .,.\right]  _{i}$ stands for the equivalence
classes of $TM_{i}$ with respect to the classical equivalence relations
between paths%
\[
\alpha\sim_{x}\beta\quad\Leftrightarrow\quad\left\{
\begin{array}
[c]{c}%
\alpha\left(  0\right)  =\beta\left(  0\right)  =x\\
\alpha^{\prime}\left(  0\right)  =\beta^{\prime}\left(  0\right)
\end{array}
\right.
\]
where $\alpha^{\prime}\left(  t\right)  =\left[  d\alpha\left(  t\right)
\right]  \left(  1\right)  $. Clearly $\Phi_{ii}=\operatorname{Id}_{TM_{i}}$
and $\Phi_{jk}\circ\Phi_{ij}=\Phi_{ik}$, for $i\leq j\leq k$, and $\{TM_{i}\}$
is a sequence of $C^{p-1}$ finite dimensional paracompact manifolds. Moreover
$\Phi_{ij}(TM_{i})$ is diffeomorphic to a closed submanifold of $TM_{j}$.

\begin{proposition}
\label{P_StructureConvenientManifold}Let $p\geq4$ and $\{M_{i},\phi_{ij}\}_{i\leq j}$ be a direct sequence of
$C^{p}$ paracompact finite dimensional manifolds for which the connecting
morphisms are $C^{p}$ embeddings with closed image.\\
Then $\varinjlim TM_{i}$ is a $C^{p-3}$ manifold modelled on $\mathbb{R}^{\infty}\times\mathbb{R}^{\infty}%
=\varinjlim(\mathbb{R}^{i}\times\mathbb{R}^{i})$.
\end{proposition}

\textbf{Proof.---}
Let $[f,x]\in\varinjlim TM_{i}$. Then for some $n\in\mathbb{N}$,
$[f,x]=\phi_{n}([f_{n},x_{n}])\in TM_{n}$. Without loss of generality suppose
that $TM_{1}\subseteq TM_{2}\subseteq\dots\subseteq TM$ and $[f,x]\in
TM_{n(x)}$. This means that $x$ belongs to $M_{n}$ and $f:(-\epsilon
,\epsilon)\longrightarrow M_{n(x)}$ is a smooth curve passing trough $x$.
Since $\{M_{i},\phi_{ij}\}_{i\leq j}$ is a directed system of manifolds
satisfying Lemma 4.1. of \cite{Glo1}, then there exists an open neighbourhood
$O_{x}$ of $x$ in $M$ and a family of $C^{p-2}$ diffeomorphisms $\{h_{i}%
^{(x)}:\mathbb{R}^{d_{i}}\rightarrow U_{i}\}_{i\geq n(x)}$ where $U_{i}%
={\phi_{i}}^{-1}(O_{x})$ and (\ref{compatibility condition of hi}) holds true.
For $i\geq n(x)$ define
\begin{align*}
Th_{i}^{(x)}:\mathbb{R}^{d_{i}}\times\mathbb{R}^{d_{i}}  &  \longrightarrow
TU_{i}\subseteq TM_{i}\\
(\bar{y},\bar{v})  &  \longmapsto\lbrack\gamma,y]
\end{align*}
where $({h_{i}^{(x)}}^{-1}\circ\gamma)(t)=\bar{y}+t\bar{v}$. For $i\leq j$ we
get
\[
\Phi_{ij}\circ Th_{i}^{(x)}(\bar{y},\bar{v})=\Phi_{ij}([\gamma,y])=[\phi
_{ij}\circ\gamma,\phi_{ij}(y)].
\]
On the other hand,
\[
Th_{j}^{(x)}\circ\left(  \lambda_{ij}\times\lambda_{ij}\right)  (y,v)=Th_{j}%
^{(x)}\left(  \left(  \bar{y},0\right)  ,\left(  \bar{v},0\right)  \right)
=[\gamma^{\prime},y^{\prime}]
\]
for which ${(h_{j}^{(x)}}^{-1}\circ\gamma^{\prime})(t)=(y,0)+t(v,0)=\lambda
_{ij}(\bar{y}+t\bar{v})$. We claim that $[\phi_{ij}\circ\gamma,\phi
_{ij}(y)]=[\gamma^{\prime},y^{\prime}]$.

Using ({\ref{compatibility condition of hi}}) we observe that
\begin{align*}
{h_{j}^{(x)}}^{-1}\circ(\phi_{ij}\circ\gamma(t))  &  =({h_{j}^{(x)}}^{-1}%
\circ\phi_{ij})\circ\gamma(t)=(\lambda_{ij}\circ{h_{i}^{(x)}}^{-1})\circ
\gamma(t)\\
&  =\lambda_{ij}\circ({h_{i}^{(x)}}^{-1}\circ\gamma)(t)=\lambda_{ij}(\bar
{y}+t\bar{v}),
\end{align*}
which proves the assertion.

Roughly speaking for any $[f,x]\in TM$, we constructed a family of $C^{p-3}$
diffeomorphisms
\[
\{Th_{i}^{(x)}:\mathbb{R}^{d_{i}}\times\mathbb{R}^{d_{i}}\longrightarrow
TU_{i}\subseteq TM_{i}\}_{i\geq n(x)}%
\]
which satisfy the compatibility conditions
\[
\Phi_{ij}\circ h_{i}^{(x)}=h_{j}^{(x)}\circ\left(  \lambda_{ij}\times
\lambda_{ij}\right)  ~;~~j\geq i\geq n(x).
\]
As a consequence the limit map $Th^{(x)}=\varinjlim Th_{i}^{(x)}%
:\mathbb{R}^{\infty}\times\mathbb{R}^{\infty}\longrightarrow TU^{(x)}:=\bigcup\limits_{i\geq n(x)}
TU_{i}$ can be defined. The map $Th^{(x)}$ denotes the diffeomorphism whose
restriction to $\mathbb{R}^{d_{i}}\times\mathbb{R}^{d_{i}}$ is $Th_{i}^{(x)}$.

The next step is to establish that the family ${\mathcal{B}}=\{{Th^{(x)}}%
^{-1};x\in M\}$ is an atlas for $TM$. For $[f,x]$ and $[f^{\prime},x^{\prime
}]$ in $TM$ define $n=\max\{n(x),n(x^{\prime})\}$. Set $\tau:=Th^{(x^{\prime
})}\circ{Th^{(x)}}^{-1}$. Since for $i\geq n$
\[
\tau\circ\lambda_{i}=\lambda_{i}\circ Th_{i}^{(x^{\prime})}\circ{Th_{i}^{(x)}%
}^{-1}%
\]
it follows that $\tau$ is a $C^{p-3}$ diffeomorphism too. Moreover for every
natural number $i$, $TM_{i}$ is a locally compact topological space. This last
means that $\varinjlim TM_{i}$ is Hausdorff (\cite{Han}, \cite{Glo2}) which
completes the proof.  \hfill $\blacksquare$

\subsection{The Lie group $Gl\left(  \infty,\mathbb{R}\right)  $}

In the situation described in \cite{Gal} (tangent bundle of projective limit
of Banach manifolds), the general linear group $GL\left(  \mathbb{F}\right)  $
cannot play the r\^{o}le of structural group and is replaced by $H_{0}\left(
\mathbb{F}\right)  $ which is a projective limit of Banach Lie groups.

In our framework we are going to use the convenient Lie group $GL\left(
\infty,\mathbb{R}\right)  $ as structural group. It is defined as follows. The
canonical embeddings $\mathbb{R}^{n}\longrightarrow\mathbb{R}^{n+1}$ induce
injections $GL\left(  \mathbb{R}^{n}\right)  \longrightarrow GL\left(
\mathbb{R}^{n+1}\right)  .$ The inductive limit is given by%
\[
GL\left(  \infty,\mathbb{R}\right)  =\varinjlim GL\left(  \mathbb{R}%
^{n}\right)
\]

and can be endowed with a real analytic regular Lie group modeled on
$\mathbb{R}^{\infty}$ (cf \cite{KriMic}, Theorem 47.8).

\subsection{Convenient vector bundle structure on $TM$}

\begin{theorem}
\label{T_StructureConvenientVectorBundleOnTM}$TM$ over $M$ admits a convenient
vector bundle structure with the structure group $GL\left(  \infty
,\mathbb{R}\right)  $.
\end{theorem}

\textbf{Proof.---}
For any $i\in\mathbb{N}$ consider the natural projection $\pi_{i}%
:TM_{i}\longrightarrow M_{i}$ which maps $[\gamma,y]$ onto $y$. As a first
step we show that the limit map $\pi:=\varinjlim\pi_{i}$ exists. For $j\geq i$
and $[\gamma,x]\in TM_{i}$ we have
\[
\phi_{ij}\circ\pi_{i}[\gamma,y]=\phi_{ij}\left(  y\right)
\]
On the other hand
\[
\pi_{j}\circ\Phi_{ij}[\gamma,y]=\pi_{j}[\phi_{ij}\circ\gamma,\phi_{ij}\left(
y\right)  ]
\]
The compatibility condition $\phi_{ij}\circ\pi_{i}=\pi_{j}\circ\Phi_{ij}$
leads us to the limit (differentiable) map
\[
\pi:=\varinjlim\pi_{i}:\varinjlim TM_{i}\longrightarrow\varinjlim M_{i}%
\]
whose restriction to $TM_{i}$ is given by $\phi_{i}\circ\pi_{i}=\pi\circ
\Phi_{i}$.

For $[f,x]\in\varinjlim TM_{i}$ consider the family of diffeomorphisms
$\{h_{i}^{(x)}:\mathbb{R}^{d_{i}}\longrightarrow U_{i}^{(x)}\}_{i\geq n(x)}$
as before. For any $i\geq n(x)$ define
\begin{align*}
\Psi_{i}:{\pi_{i}}^{-1}(U_{i}^{(x)})  &  \longrightarrow U_{i}^{(x)}%
\times\mathbb{R}^{d_{i}}\\
{[\gamma,y]}  &  \longmapsto\left(  y,({h_{i}^{(x)}}^{-1}\circ\gamma)^{\prime
}(0)\right)  .
\end{align*}
With the standard calculation for the finite dimensional manifolds it is known
that $\Psi_{i}$, $i\in\mathbb{N}$, is a diffeomorphism. For $j\geq i\geq
n(x)$, we claim that the following diagram is commutative
\[%
\begin{array}
[c]{ccc}%
{\pi_{i}}^{-1}(U_{i}^{(x)}) & \overset{\Psi_{i}}{\rightarrow} & U_{i}%
^{(x)}\times\mathbb{R}^{d_{i}}\\
\Phi_{ij}\downarrow &  & \downarrow\phi_{ij}\times\lambda_{ij}\\
{\pi_{j}}^{-1}(U_{j}^{(x)}) & \overset{\Psi_{j}}{\rightarrow} & U_{j}%
^{(x)}\times\mathbb{R}^{d_{j}}%
\end{array}
\]

To see that we argue as follows.
\begin{align*}
(\phi_{ij}\times\lambda_{ij})\circ\Psi_{i}([\gamma,y])  &  =(\phi_{ij}%
\times\lambda_{ij})\left(  y,({h_{i}^{(x)}}^{-1}\circ\gamma)^{\prime
}(0)\right) \\
&  =\left(  \phi_{ij}(y),\lambda_{ij}\circ\left(  ({h_{i}^{(x)}}^{-1}%
\circ\gamma)^{\prime}(0)\right)  \right) \\
&  \overset{\left(  \ast\right)  }{=}\left(  \phi_{ij}(y),(\lambda_{ij}%
\circ{h_{i}^{(x)}}^{-1}\circ\gamma)^{\prime}(0)\right) \\
&  \overset{\left(  \ast\ast\right)  }{=}\left(  \phi_{ij}(y),\left(
({h_{j}^{(x)}}^{-1}\circ\phi_{ij}\circ\gamma)^{\prime}(0)\right)  \right) \\
&  =\Psi_{j}\left(  [\phi_{ij}\circ\gamma,\phi_{ij}(y)]\right) \\
&  =\left(  \Psi_{j}\circ\Phi_{ij}\right)  \left[  \gamma,y\right]
\end{align*}

For $(\ast\ast)$ we used the equation (\ref{compatibility condition of hi})
and for $(\ast)$ using the linearity of $\lambda_{ij}$\\
we get
\begin{align*}
\lambda_{ij}\circ\left(  ({h_{i}^{(x)}}^{-1}\circ\gamma)^{\prime}(0)\right)
&  =\lambda_{ij}\left(  \underset{t\rightarrow0}{\lim}\dfrac{({h_{i}^{(x)}%
}^{-1}\circ\gamma)(t)-({h_{i}^{(x)}}^{-1}\circ\gamma)(0)}{t}\right) \\
&  =\underset{t\rightarrow0}{\lim}\dfrac{(\lambda_{ij}\circ{h_{i}^{(x)}}%
^{-1}\circ\gamma)(t)-(\lambda_{ij}\circ{h_{i}^{(x)}}^{-1}\circ\gamma)(0)}{t}\\
&  =(\lambda_{ij}\circ{h_{i}^{(x)}}^{-1}\circ\gamma)^{\prime}(0).
\end{align*}

Since $\pi_{i}^{-1}(U_{i}^{(x)})$, $i\geq n(x)$, is open and since $\pi
^{-1}(U)=\varinjlim\pi_{i}^{-1}(U_{i}^{(x)})$, it follows that $\pi
^{-1}(U)\subseteq TM$ is open. Furthermore $\Psi_{x}:=\varinjlim\Psi_{i}%
:\pi^{-1}(U)\longrightarrow U\times\mathbb{R}^{\infty}$ exists and, as a
direct limit of $C^{p-3}$ diffeomorphisms, is a $C^{p-3}$ diffeomorphism. On
the other hand
\[
\Psi_{x}|_{\pi^{-1}(y)}:{\pi^{-1}(y)}\longrightarrow\{y\}\times\mathbb{R}%
^{\infty}%
\]
is linear and $pr_{1}\circ\Psi_{x}$ coincides with $\pi$. ($pr_{1}$ stands for
projection to the first factor.)

Suppose that $[f,x],[g,y]\in TM$, $n=\max\{n(x),n(y)\}$ and the intersection
$U_{xy}:=U^{(x)}\cap U^{(y)}$ is not empty. Then
\[
\left(  \Psi_{y}\right)  ^{-1}|_{U_{xy}\times\mathbb{R}^{\infty}}\circ\Psi
_{x}|_{U_{xy}\times\mathbb{R}^{\infty}}:{U_{xy}\times\mathbb{R}^{\infty}%
}\longrightarrow{U_{xy}\times\mathbb{R}^{\infty}}%
\]

arises as the inductive limit of the family
\begin{align*}
\left(  {\Psi_{i}^{y}}\right)  ^{-1}|_{U_{i}^{xy}\times\mathbb{R}^{d_{i}}%
}\circ\Psi_{i}^{x}|_{U_{i}^{xy}\times\mathbb{R}^{d_{i}}}:U_{i}^{xy}  &
\longrightarrow GL\left(  \mathbb{R}^{d_{i}}\right) \\
\bar{y}  &  \longmapsto T_{xy}^{i}(\bar{y}).
\end{align*}
Finally the family of maps $\{T_{xy}^{i}:=\left(  {\Psi_{i}^{y}}\right)
^{-1}|_{U_{i}^{xy}\times\mathbb{R}^{d_{i}}}\circ\Psi_{i}^{x}|_{U_{i}%
^{xy}\times\mathbb{R}^{d_{i}}}\}_{i\geq n}$, satisfy the required
compatibility condition and their limit $T_{xy}:=\varinjlim T_{xy}^{i}$
belongs to $\varinjlim GL\left(  \mathbb{R}^{d_{i}}\right)  :=GL\left(
\infty,\mathbb{R}\right)  $. \\
Consequently $\varinjlim TM_{i}$ becomes a(convenient) vector bundle with the fibres of type $\mathbb{R}^{\infty}$ and
the structure group $GL\left(  \infty,\mathbb{R}\right)  $. \hfill $\blacksquare$

\begin{example}
Tangent bundle to $\mathbb{S}^{\infty}$.-- The tangent bundle $T\mathbb{S}%
^{\infty}$ to the sphere $\mathbb{S}^{\infty}$ is diffeomorphic to $\left\{
\left(  x,v\right)  \in\mathbb{S}^{\infty}\times\mathbb{R}^{\infty
}:\left\langle x,v\right\rangle =0\right\}  $.
\end{example}

\begin{proposition}
\label{P_IsomorphismDirectLimitTangentBunddle}
$\varinjlim TM_{i}$ as a set is isomorphic to $T M$.
\end{proposition}

\textbf{Proof.---}
Arguing as before, let $[f,x]\in\varinjlim TM_{i}$. Then there exists
$n(x)\in\mathbb{N}$ such that, for $i\geq n(x)$, $[f,x]$ belongs to $TM_{i}$
which means that $x\in M_{i}$ and $f:(-\epsilon,\epsilon)\longrightarrow
M_{i}$ for some $\epsilon>0$. This last means that $f:(-\epsilon
,\epsilon)\longrightarrow\varinjlim M_{i}$ and consequently $[f,x]$ belongs to
$T M$.

Conversely, suppose that $[f,x]$ belongs to $TM$ that is $x\in M$ and $f$ is a
curve in $M=\varinjlim M_{i}$. Again there exists $n(x)$ such that $x\in
M_{i}$ and $f:(-\epsilon,\epsilon)\longrightarrow M_{i}$ is a smooth curve for
$i\geq n(x)$. Since $[f,x]\in TM_{i}$, $i\geq n(x)$, then $[f,x]\in\varinjlim
TM_{i}$ which completes the proof.
\hfill $\blacksquare$

\newpage

Ali Suri \\
Department of Mathematics, Faculty of Science, Bu Ali Sina University, Hamedan, 65178, Iran\\
E-mail: a.suri@math.iut.ac.ir\\

Patrick Cabau  \\
Lyc\'{e}e Pierre de Fermat, 2 Parvis des Jacobins, BP 7013,  31068 Toulouse Cedex 7, France\\
E-mail: Patrick.Cabau@ac-toulouse.fr\\


\begin{thebibliography}{KriMic}

\bibitem[ADGS]{ADGS} M. Aghasi, C.T. Dodson, G.N. Galanis, A. Suri, \textit{%
Conjugate connections and differential equations on infinie dimensional
manifolds}, J. Geom. Phys. (2008)

\bibitem[AghSur]{AghSur} M. Aghasi, A. Suri, \textit{Ordinary differential
equations}, Balkan Journal of Geometry and Its Applications \textbf{12} n${%
{}^{\circ }}2$ (2007) 1--8

\bibitem[Bou]{Bou}N. Bourbaki, \textit{El\'{e}ments de Math\'{e}matiques},
Alg\`{e}bre, Chapitres 1 \`{a} 3, 2$^{\grave{e}me}$ \'{e}dition, Springer 2006

\bibitem[Gal]{Gal}G.N. Galanis, \textit{Differential and Geometric Structure
for the Tangent Bundle of a Projective Limit Manifold}, Rend. Sem. Univ.
Padova, Vol. 112 (2004)

\bibitem[Glo1]{Glo1}H. Gl\"{o}ckner, \textit{Direct limit of Lie groups and
manifolds}, J. Math. Kyoto Univ. (JMKYAZ) 43-1 (2003), 1--26

\bibitem[Glo2]{Glo2}H. Gl\"{o}ckner, \textit{Fundamentals of Direct Lie
Theory}, Compositio Math. 141 (2005), 1551--1577

\bibitem[Han]{Han}V.L. Hansen, \textit{Some Theorems on Direct Limit of
Expanding Sequences of Manifolds}, Math.\ Scand.29 (1971), 5--36

\bibitem[KriMic]{KriMic}A. Kriegel, P.W. Michor, \textit{The convenient
Setting of Global Analysis }(AMS Mathematical Surveys and Monographs)
\textbf{53} 1997

\bibitem[Lan]{Lan}S. Lang, \textit{Differential and Riemannian Manifolds},
Graduate Texts in Mathematics, 160, Springer, New York 1995

\bibitem[Mil]{Mil}J. Milnor, \textit{On infinite dimensional Lie groups}, Preprint, Institute of Advanced Study, Princeton, 1982

\end{thebibliography}
\end{document}